\documentclass[a4paper,12pt]{article}
\usepackage{amsmath,amssymb,amsthm}
\usepackage{natbib}

\makeatletter
\newcommand{\Macro}[2]{\@namedef{clm@macro@#1}{#2}}
\newcommand{\m}[1]{\@nameuse{clm@macro@#1}}

\newcommand{\equal}{\@ifstar{\clm@starred@equal}{\clm@nonstarred@equal}}%
\newcommand{\clm@starred@equal}[1]{\stackrel{#1}{=}}%
\newcommand{\clm@nonstarred@equal}[1]{\stackrel{\m{#1}}{=}}

\global\def\clm@size@auto{auto}
\newcommand{\clm@casesize}[4]{\def\clm@size{#1}%
\ifx\clm@size\clm@size@auto #2\else\@ifempty{#1}{#3}{#4}\fi}

\@ifundefined{minaw@}{\newdimen\minaw@\minaw@11pt}{}

\newcommand{\xrightglobalarrow}[4]{%
  \mathrel{\mathop{%
    \setbox\z@\vbox{\m@th
      \hbox{$\scriptstyle\;{#4}\;\;$}%
      \hbox{$\m@th\scriptstyle\;{#3}\;\;$}}
    \hbox to\ifdim\wd\z@>\minaw@\wd\z@\else\minaw@\fi{%
\m@th\setboxz@h{$\displaystyle#1$}\ht\z@\z@
  $\displaystyle\copy\z@\mkern-6mu\cleaders
  \hbox{$\displaystyle\mkern-2mu\box\z@\mkern-2mu$}\hfill
  \mkern-6mu\mathord#2$}}
  \limits^{#3}\@ifnotempty{#4}{_{#4}}}}
\newcommand{\Conv}{\@ifstar{\clm@starred@Conv}{\clm@nonstarred@Conv}}
\newcommand{\clm@starred@Conv}[2]{\xrightglobalarrow{\Relbar}{\Rightarrow}{#1}{#2}}
\newcommand{\clm@nonstarred@Conv}[2]{\xrightglobalarrow{\Relbar}{\Rightarrow}{\m{#1}}{#2}}
\newcommand{\twoarg}[7][auto]{%
\def\clm@vert@bar{|}%
\edef\clm@arg@bar{#7}%
\ifx\clm@arg@bar\clm@vert@bar%
\clm@casesize{#1}{\left#4#6 #2\;\vrule\; #3#6\right#5}{%
#4#6 #2\;\vert\;#3#6 #5}{%
\@nameuse{#1l}#4#6 #2\;\@nameuse{#1}\vert\;#3#6 \@nameuse{#1r}#5}
\else
\clm@casesize{#1}{\left#4#6 #2 #7 #3 #6\right#5}{%
#4#6 #2 #7 #3 #6 #5}{%
\@nameuse{#1l}#4 #6 #2 #7 #3 #6\@nameuse{#1r}#5}
\fi}
\newcommand{\onearg}[5][auto]{%
\clm@casesize{#1}{\left#3#5 #2#5 \right#4}{%
#3#5 #2#5 #4}{\@nameuse{#1l}#3#5 #2#5 \@nameuse{#1r}#4}}

\newcommand{\clm@norm}[3]{%
\clm@casesize{#1}{\left\lVert #3\right\rVert\@ifnotempty{#2}{_{#2}}}{%
\lVert #3\rVert\@ifnotempty{#2}{_{#2}}}{%
\@nameuse{#1l}\lVert #3\@nameuse{#1r}\rVert\@ifnotempty{#2}{_{#2}}}}
\newcommand{\clm@starred@norm}[3][auto]{\clm@norm{#1}{#2}{#3}}
\newcommand{\clm@nonstarred@norm}[3][auto]{\@ifundefined{clm@named@norm@#2}%
{\clm@norm{#1}{\m{#2}}{#3}}{\@nameuse{clm@named@norm@#2}{#1}{#3}}}
\newcommand{\norm}{\@ifstar{\clm@starred@norm}{\clm@nonstarred@norm}}

\newcommand{\scalprod}{\@ifstar{\clm@starred@scalprod}%
{\clm@nonstarred@scalprod}}
\newcommand{\clm@starred@scalprod}[4][auto]{
\twoarg[#1]{#3}{#4}{\langle}%
{\rangle}{}{,}\@ifempty{#2}{}{_{#2}}}
\newcommand{\clm@nonstarred@scalprod}[4][auto]{
\@ifundefined{clm@scalprod@named@#2}{
\twoarg[#1]{#3}{#4}{\langle}%
{\rangle}{}{,}\@ifempty{#2}{}{_{\m{#2}}}}{\@nameuse{clm@scalprod@named@#2}{#1}{#3}{#4}}}
\makeatother

\newcommand{\ens}[3][auto]{\twoarg[#1]{#2}{#3}\{\}{\,}|}
\newcommand{\event}[2][auto]{\onearg[#1]{#2}\{\}{\,}}

\Macro{law}{\mathrm{dist.}}
\Macro{proba}{\mathrm{prob.}}

\DeclareMathOperator{\Dom}{Dom}
\newcommand{\cB}{\mathcal{B}}
\newcommand{\cU}{\mathcal{U}}
\newcommand{\cV}{\mathcal{V}}
\newcommand{\cC}{\mathcal{C}}
\newcommand{\cX}{\mathcal{X}}
\newcommand{\cW}{\mathcal{W}}
\newcommand{\fS}{\mathfrak{S}}

\newcommand{\hp}{\hat p}
\newcommand{\hq}{\hat q}

\newcommand{\hwn}{\hat W^n}
\newcommand{\ot}{[0,T]}
\newcommand{\ou}{[0,1]}
\newcommand{\snu}{\sum_{n\ge 1}}
\newcommand{\wo}{\mathcal{W}_0}

\Macro{Id}{\mathrm{Id}}
\Macro{cBalpha}{\cB_\alpha}
\Macro{cBrho}{\cB_\rho}
\Macro{cU->Brho}{\cU\to\cB_\rho}
\Macro{cBrho'}{\cB'_\rho}
\Macro{cBdelta}{\cB_\delta}
\Macro{cB-alpha}{\cB_{-\alpha}}
\DeclareMathOperator{\rg}{rg}
\Macro{Lip}{\mathrm{Lip}}
\Macro{cV}{\cV}
\Macro{cX}{\cX}

\newcommand{\RR}{\mathbb{R}}
\newcommand{\NN}{\mathbb{N}}
\newcommand{\vd}{\,\mathrm{d}}
\newcommand{\dd}{\mathrm{d}}

\newcommand{\cb}{\mathcal B}

\newcommand{\cf}{\mathcal F}

\newcommand{\cl}{\mathcal L}

\newcommand{\al}{\alpha}
\newcommand{\ep}{\varepsilon}

\newcommand{\si}{\sigma}

\newcommand{\vp}{\varphi}

\newcommand{\lp}{\left(}
\newcommand{\rp}{\right)}
\newcommand{\lc}{\left[}
\newcommand{\rc}{\right]}

\newtheorem{theorem}{Theorem}
\newtheorem{proposition}{Proposition}
\newtheorem{lemma}{Lemma}

\newtheorem{hypothesis}{Hypothesis}
\theoremstyle{definition}

\theoremstyle{remark}
\newtheorem{remark}{Remark}

\begin{document}

\title{Young integrals and SPDEs}
\author{
  \\
  { Antoine Lejay}              \\
  {\small\it Project OMEGA, INRIA Lorraine }  \\[-0.2cm]
  {\small\it IECN, Campus Scientifique}          \\[-0.2cm]
  {\small\it BP 239 -- 54506 Vand\oe uvre-l\`es-Nancy CEDEX, France}  \\[-0.2cm]
  {\small  {\tt Antoine.Lejay@iecn.u-nancy.fr}}   
  \\[-0.1cm]  
\and
  { Massimiliano Gubinelli}              \\
  {\small\it  Laboratoire Analyse, G\'eom\'etrie \& Applications --  UMR 7539 
}  \\[-0.2cm]
  {\small\it  Institute Galil\'ee, Universit\`e Paris 13,}  \\[-0.2cm]
  {\small\it 93430 -- Villetaneuse, France}      \\[-0.2cm]
  {\small  {\tt gubinell@math.univ-paris13.fr}}\\[-0.1cm]  
\and
  { Samy Tindel}              \\
  {\small\it Institut Elie Cartan} \\[-0.2cm]
  {\small\it  Universit\'e Henri Poincar\'e (Nancy)}          \\[-0.2cm]
  {\small\it BP 239 -- 54506 Vand\oe uvre-l\`es-Nancy CEDEX, France}  \\[-0.2cm]
  {\small  {\tt tindel@iecn.u-nancy.fr }} \\[-0.1cm]  
  {\protect\makebox[5in]{\quad}}  
  \\
}


\maketitle

\begin{abstract}
In this note, we study the non-linear evolution problem
$$
dY_t = -A Y_t dt + B(Y_t) dX_t,
$$  
where $X$ is a $\gamma$-H\"older continuous function of the time parameter, with values in a distribution space, and $-A$ the
generator of an analytical semigroup. Then, we will give some sharp conditions on $X$ in order to solve the above equation in a function space, first in the linear case (for  any value of $\gamma$ in $(0,1)$), and then when $B$ satisfies some Lipschitz type conditions (for $\gamma>1/2$). The solution of the evolution problem will be understood in the {\it mild} sense, and the integrals involved in that definition will be of Young type.
\end{abstract}
\section{Introduction}
The existence and uniqueness problem for ordinary differential equations driven by an irregular path of H\"older continuity greater than $1/2$ is now fairly well understood, either in the context of fractional integrals \cite{Za}, or as a first step towards the definition of differential equations driven by a rough path (see \cite{Gubinelli,Lej,LyonsBook}).

\vspace{0.3cm}

However, the case of partial differential equations of evolution type has only been partially treated. For instance, when the driving noise is an infinite dimensional fractional Brownian motion of Hurst parameter $H$, let us mention two results that have been obtained recently:
\begin{itemize}
\item
In case of a linear equation with additive noise, some optimal conditions  on the space covariance of the noise, ensuring the existence and uniqueness of a function-valued solution to the evolution equation, are given in \cite{TTV}. These results are based on a representation of the solution involving some Skorokhod type integrals, that cannot be interpreted as the limit of Riemann sums.
\item
In the non-linear case, the same kind of results are available in \cite{Maslo-Nualart}. In this latter case, only the case $H>1/2$ is considered, the space covariance is assumed to be trace class (which implies a strong regularity condition on the driving noise) and the integrals are defined path-wise.
\end{itemize}
On the other hand, some differential equations in a Banach space, driven by a rough path, are considered in \cite{Ledoux-Lyons}, but the question of the regularization of the noise by a semi-group is not addressed in that last reference.

\vspace{0.3cm}

With all those facts in mind, our aim, in this paper, is to make a step
towards the resolution of non linear partial differential equations driven by
an infinite dimensional rough path. To be more specific, let us start by
defining roughly the kind of equation we handle: we consider  an
unbounded operator $A$ on a Banach space $\cb$, and assume that $-A$ is the
infinitesimal generator of an analytical semi-group $(S(t))_{t\geq 0}$. This
induces a family of Banach spaces $\{ \cb_\al;
\al\in\mathbb{R} \}$, where, for $\al>0$, $\cb_\al=\Dom(A^{\al})$, and
$\cb_{-\al}$ is defined by duality. Now, our noise $X$ is a function from
$[0,T]$ to $\cb_{-\al}$, for a certain $\al>0$, with a given H\"older
continuity $\gamma\in(0,1)$ in time. Our equation is of the form
\begin{equation}\label{defevol}
\vd Y_t = -A Y_t \vd t + B(Y_t) \vd X_t, \quad t\in [0,T],
\end{equation}
with a given initial condition $y_0$, and where $B$ is a map from $\cb_\delta$
to $\cl(\cb_{-\al};\cb_\delta)$ satisfying some (local) Lipschitz conditions.
In fact, equation (\ref{defevol}) is understood in the so-called {\it
mild} sense, i.e. we say that $Y$ is a solution  to (\ref{defevol}) in
$[0,T]$ if it is a function in $C^{\kappa}([0,T];\cb_{\delta})$, with  a given
$\kappa>0$, satisfying 
$$
Y(t)=S(t)y_0
+\int_0^t S(t-s) B(Y(s))\, \vd X(s), \quad t\le T.
$$
Based on these notations, we get the following results:
\begin{enumerate}
\item
In the linear additive case, that is when $B=\mbox{Id}$, we get the existence
and uniqueness of a global solution to (\ref{defevol}), living in
$C^{\kappa}([0,T];\cb_{\delta})$ for any $T>0$, with the following condition on
the coefficients $\al,\delta,\gamma$ and $\kappa$:
$$
\al<\gamma,\quad
0<\delta<\gamma-\al,\quad
0<\kappa<\gamma-\al-\delta.
$$
This result is given in a rigorous form at Theorem \ref{thm-0}.
\item
In the general case, if $B$ satisfies some Lipschitz type conditions, and under
the additional assumption $\gamma+\delta>1$ (which implies in particular that
$\gamma>1/2$), we also get the existence and uniqueness of the solution to
(\ref{defevol}). In case of a locally Lipschitz coefficient $B$, the existence
of a solution can only be guaranteed up to an explosion time $T>0$.
\item
Eventually, we check that our abstract results can be applied to a simple
case, namely the case of the heat equation in $\ou$ with Dirichlet boundary
conditions. We also assume that $X$ is a cylindrical fractional Brownian
noise (see Section \ref{defrac} for a complete description). Then we are
able to solve the equation (\ref{heateq}) for a function $\si\in
C_b^2(\mathbb{R})$, up to a strictly positive explosion time. At that point,
let us insist again on the following fact: we are allowed to consider, in the
case of equation (\ref{heateq}), a white noise in space for any $H>1/2$. To our
knowledge, this is the first occurrence of an existence and uniqueness result
for a non linear SPDE driven by a fractional Brownian motion with a non trivial
space covariance.
\end{enumerate}
Let us also observe that, in order to get those results, we use a quite
natural approach: our setting allows us to use Young integrals, that can be
approximated by Riemann sums. Then, in each of those Riemann sums, 
we regularize the path $X_t\in\cb_{-\al}$ by the
semi-group $S$ in order to get
an element of $\cb_{\delta}$. This induces an additional singularity in
time, that we compensate by the H\"older regularity of~$X$. 
Of course, one may want to consider some more irregular paths, in which
case one could expect to need 
paths living in a ``bigger'' space that $\cb_{-\al}$,
as in the rough path theory. But we believe that 
the results contained in this article may represent 
one of the corner stone for the development of such a general theory.

\vspace{0.3cm}

Our paper is organized as follows: at Section \ref{addcase}, we recall
some basic facts about analytic semi-groups, and treat the linear case. Section
\ref{nonlin} is devoted  to the non-linear case, and Section \ref{exheat},
to the example of the heat equation in dimension 1.

\section{Case of an additive noise}\label{addcase}

In this Section, we first recall some basic facts about analytical semi-groups, that we will use throughout the paper, and then solve equation (\ref{defevol}) in the linear additive case.

\subsection{Analytical semi-group}

This section contains some classical results
about analytical semi-groups and fractional power
of their infinitesimal generators. 
For proofs, 
see \cite{pazy83a,engel00a,fattorini99a} for example.

Let $(\cB,\norm{}{\cdot})$ be a separable Banach space.
Let $(A,\Dom(A))$ be a non-bounded linear operator 
on $\cB$ such that $-A$ be the infinitesimal generator of
an analytical semi-group $(S(t))_{t\geq 0}$.
Assume that for some constants $M>0$ and $\lambda<0$, 
$\norm*{\cB\to\cB}{S(t)}\leq Me^{-\lambda t}$ for all $t\geq 0$.
In particular, this implies that $A$ is one-to-one from
$\Dom(A)$ to $\cB$.

For any $\alpha\in\RR$, the fractional power $(A^\alpha,\Dom(A^\alpha))$
of $A$ can be defined. If $\alpha<0$, then 
$\Dom(A^\alpha)=\cB$ and 
$A^\alpha$ is one-to-one from $\cB$ to $\rg(A^{\alpha})$
and is a bounded operator. Besides, if $\alpha\geq 0$,
$(A^\alpha,\Dom(A^\alpha))$
is a closed operator
with a dense domain $\Dom(A^\alpha)=\rg(A^{-\alpha})$.
Moreover, $A^\alpha=(A^{-\alpha})^{-1}$.

For $\alpha\geq 0$, let $\m{cBalpha}$ be the space
$\Dom(A^\alpha)$ with the norm 
$\norm{cBalpha}{x}=\norm{}{A^\alpha x}$.
Since $A^{-\alpha}$ is continuous, 
it follows that the norm $\norm{cBalpha}{\cdot}$ 
is equivalent to the graph norm of 
$A^\alpha$.
If $\alpha=0$, then $\m{cBalpha}=\cB$
and $A^0=\m{Id}$.

If $\alpha<0$, let $\m{cBalpha}$ be
the completion of $\cB$ with respect 
to $\norm{cBalpha}{x}=\norm{}{A^\alpha x}$.
Thus, $\m{cBalpha}$ is a larger space than 
$\cB$.

Among the important facts about these spaces, 
note the following ones:
For any $\alpha\in\RR$ and any $\rho\geq 0$, 
\begin{gather}
A^{-\rho}\text{ maps }\m{cBalpha}\text{ onto }
\cB_{\alpha+\rho}\text{, for all }\alpha\in\RR,\ \rho\geq 0,\\
\norm{cBalpha}{x}\leq C_{\alpha,\delta}\norm{cBdelta}{x}
\text{ for all }x\in\m{cBalpha}
\text{ and all }\alpha\leq \delta.
\end{gather}
Moreover, for all $\alpha,\beta\in\RR$,
\begin{equation}
\label{eq-sg-5}
A^\alpha A^\beta=A^{\alpha+\beta}
\text{ on }\cB_\gamma
\end{equation}
with $\gamma=\max\event{\alpha,\beta,\alpha+\beta}$.

The semi-group $(S(t))_{t\geq 0}$ also satisfies 
\begin{gather}
\label{eq-sg-1}
S(t)\text{ may be extended to $\m{cBalpha}$ for all $\alpha<0$ 
and all $t\geq 0$},\\
\label{eq-sg-2}
S(t)\text{ maps }\m{cBalpha}\text{ to }\m{cBdelta}
\text{ for all }\alpha\in\RR,\ \delta\geq 0,\\
\label{eq-sg-3}
\text{ for all $t>0$, $\alpha\in\RR$, } 
\norm*{\cB\to\cB}{A^\alpha S(t)}\leq M_\alpha t^{-\alpha}e^{-\lambda t},\\
\label{eq-sg-4}
\text{ for $0<\alpha\leq 1$, $x\in\m{cBalpha}$, }
\norm{}{S(t)x-x}\leq C_\alpha t^\alpha\norm{}{A^\alpha x}
\end{gather}
where with $\norm*{\cB\to\cB}{\cdot}$ we denote the operator norm from
$\cB$ to $\cB$.

\subsection{Mild solutions of the Cauchy problem with an additive noise}

This subsection is devoted to the linear additive case
of our evolution equation, for which we will
introduce first some additional notation:
for any function $f$ defined on $\RR_+$, set
$f(s,t)=f(t)-f(s)$. For a $\gamma$-H\"older continuous
function $f$ from $[0,T]$ to a Banach space $\cX$, define 
$H_{\gamma,T}(f;\cX)$ by
\begin{equation*}
H_{\gamma,T}(f;\cX)=\sup_{0\leq s\leq t\leq T}
\frac{\norm{cX}{f(s,t)}}{(t-s)^\gamma}.
\end{equation*}

Let $Y$ be a $\gamma$-H\"older continuous
function from $[0,T]$ to $\cB$ with $\gamma\in(0,1)$,
and set $H(Y)=H_{\gamma,T}(Y;\cB)$.
Fix also $\alpha\in\RR$
and set $X=A^\alpha Y$. Note that $X$ 
belongs to $\m{cB-alpha}$, but does not 
necessarily belong to $\cB$.

Consider now the following Cauchy linear equation
\begin{equation}
\label{eq-cle}
u(t)=X(t)-\int_0^t Au(s)\vd s.
\end{equation}
Without loss of generality, it can be assumed
that $X(0)=0$. Otherwise, the solution $u$
to \eqref{eq-cle} is the sum of $v$
and $w$ that are solutions to 
$v(t)=X(t)-X(0)-\int_0^t Av(s)\vd s$
and $w(t)=X(0)-\int_0^t Aw(s)\vd s$.
The solution of the later problem 
is $w(t)=S(t)X(0)$. If $X$ is smooth enough
and $X'$ belongs to $\mathrm{L}^1([0,T];\cB)$, 
then
\begin{equation*}
v(t)=\int_0^t S(t-s)X'(s)\vd s
\end{equation*}
belongs to $\cC([0,T];\cB)$ 
and is called a \emph{mild solution} of \eqref{eq-cle}.
This notion of solution is weaker than 
the notion of \emph{strong solutions}, since
$v$ is not necessarily differentiable. On
the relation between strong and mild solutions, 
see for example \cite{pazy83a}.

We can now state the main result of this section:
\begin{theorem}\label{thm-0}
Assume
that $\alpha<\gamma$.
There exists a linear map
\begin{equation*}
\fS:
\cC^\gamma([0,T];\m{cB-alpha})
\to \cC^\kappa([0,T];\m{cBdelta})
\end{equation*}
for all $\delta\in(0,\gamma-\alpha)$
and  all $\kappa\in(0,(\gamma-\alpha-\delta)\wedge1)$ 
such that, if $X$ belongs to $\cC^{1}([0,T];\m{cB-alpha})$,
\begin{equation*}
\fS(X)(t)=\int_0^t S(t-s)X'(s)\vd s.
\end{equation*}
Moreover, for all $T>0$, there exist some  constants $C_1$ and $C_2$
depending only on $\alpha$, $\kappa$, $\delta$ and $\gamma$
such that 
\begin{gather}
\label{eq-8}
H_{\kappa,T}(\fS(X);\m{cBdelta})\leq C_1
H_{\gamma,T}(X;\m{cB-alpha})\\
\label{eq-9}
\text{ and }
\sup_{t\in[0,T]}\norm*{\m{cBdelta}}{\fS(X)(t)}
\leq C_2T^{\gamma-\delta-\alpha}
H_{\gamma,T}(X;\m{cB-alpha}).
\end{gather}
\end{theorem}

\begin{proof}

Fix  $T\geq t\geq s'>s\geq 0$. Let $n\in\NN$ and set 
\begin{equation*}
\fS^n(X)(s,s';t)=\sum_{k={\lceil 2^ns/t\rceil}}^{{\lfloor 2^n s'/t\rfloor}-1} S(t-t_k^n)X(t^n_k,t^n_{k+1}),
\end{equation*}
where $t_k^n=tk/2^n$ and $X(a,b) = X(b) - X(a)$.
Using the semi-group property, 
the difference between $\fS^n(X)(s,t)$ and $\fS^n(X)(s,t)$
is 
\begin{multline*}
\fS^{n+1}(X)(s,s';t)-\fS^n(X)(s,s';t)\\
=\sum_{k={\lceil 2^ns/t\rceil}}^{{\lfloor 2^n s'/t\rfloor}-1}
(S(t-t_{2k+1}^{n+1})-S(t-t_{2k}^{n+1}))
X(t_{2k+1}^{n+1},t_{2k+2}^{n+1})\\
=\sum_{k={\lceil 2^ns/t\rceil}}^{{\lfloor 2^n s'/t\rfloor}-1}
(\m{Id}-S(t_{2k+1}^{n+1}-t_{2k}^{n+1}))
S(t-t_{2k+1}^{n+1})X(t_{2k+1}^{n+1},t_{2k+2}^{n+1})\\
\end{multline*}

Now, fix $\delta<\gamma-\alpha$ and $\beta\in (1-\gamma,(1-\alpha-\delta)\wedge 1)$.
With \eqref{eq-sg-5} and \eqref{eq-sg-3}--\eqref{eq-sg-4}, 
\begin{multline*}
\norm{}{A^\delta \fS^{n+1}(X)(s,s';t)-A^\delta\fS^n(X)(s,s';t)}\\
\leq \sum_{k={\lceil 2^ns/t\rceil}}^{{\lfloor 2^n s'/t\rfloor}-1}
C_\beta \left(\frac{t}{2^{n+1}}\right)^\beta
\norm{}{A^{\beta+\delta+\alpha}S(t-t_{2k+1}^{n+1})Y(t_{2k+1}^{n+1},t_{2k+2}^{n+1})}\\
\leq \sum_{k={\lceil 2^ns/t\rceil}}^{{\lfloor 2^n s'/t\rfloor}-1}
C_\beta M_{\alpha+\beta+\gamma}
\left(\frac{t}{2^{n+1}}\right)^{\beta+\gamma}
\frac{H(Y)}{(t-t_{2k+1}^{n+1})^{\beta+\delta+\alpha}}.
\end{multline*}
Set $\epsilon=\alpha+\delta+\beta$. By assumption, $\epsilon<1$,
and furthermore, 
\begin{align*}
\sum_{k={\lceil 2^ns/t\rceil}}^{{\lfloor 2^n s'/t\rfloor}-1}
\frac{1}{(t-t^{n+1}_{2k+1})^\epsilon}
&\leq \frac{2^{n+1}}{t^\epsilon}
 \int_{{\lceil 2^ns/t\rceil}/2^{n+1}}^{{\lfloor 2^n s'/t\rfloor}/2^{n+1}}\frac{\vd r}{(1-r)^\epsilon}\\
&\leq \frac{2^{n+1}}{t^\epsilon}
\int_0^1 \frac{\vd r}{(1-r)^{\epsilon}}.
\end{align*}
Hence, there exists a constant $C$ that depends only on $\alpha$, $\beta$ 
and $\delta$ for which 
\begin{equation}
\label{eq-4}
\norm{}{A^\delta \fS^{n+1}(X)(s,s';t)-A^\delta\fS^n(X)(s,s';t)}
\leq \frac{CH(Y)t^{\gamma-\delta-\alpha}}{2^{n(\beta+\gamma-1)}}
\end{equation}
for all integer $n$. Since $\beta+\gamma>1$, the series 
\begin{equation*}
\sum_{n\geq 0}\norm{}{A^\delta \fS^{n+1}(X)(s,s';t)-A^\delta\fS^n(X)(s,s';t)}
\end{equation*}
is convergent, and moreover, 
\begin{equation}
\label{eq-7}
\norm[]{}{A^\delta\fS^0(X)(s,s';t)}
=\norm[]{}{A^{\delta+\alpha}S(t-s)Y(s,s')}
\leq \frac{(s'-s)^\gamma H(Y)M_{\alpha+\delta}}{(t-s)^{\delta+\alpha}}.
\end{equation}
Thus, the sequence $(A^\delta\fS^n(X)(s,s';t))_{n\in\NN}$
is convergent in $\cB$. Let now $\fS(X)(s,s';t)$ be the
limit of $(\fS^n(X)(s,s';t))_{n\in\NN}$. 
Observe that, since $A^\delta$
is a closed operator, 
the limit of 
$(A^\delta\fS^n(X)(s,s';t))_{n\in\NN}$ is 
$A^\delta\fS(X)(s,s';t)$.

To simplify the notations, 
set $\fS(X)(t)=\fS(X)(0,t;t)$.
The linearity of $X\mapsto \fS(X)(t)$ 
follows immediately from the construction
of $\fS(X)(t)$. The inequality \eqref{eq-9} 
is easily obtained from \eqref{eq-4} and \eqref{eq-7}.

Note also that if $X$ is smooth,
then $\fS(X)(s,s';t)=\int_s^{s'} S(t-s)X'(s)\vd s$
for all $0\leq s\leq s'\leq t$.
Moreover  $t\mapsto \fS(X)(s,s';t)$
and $t\mapsto \fS(X)(s,t;t)$ are continuous.
If $X\in\cC^\gamma([0,T];\m{cB-alpha})$,
then one can find a sequence $(X^n)_{n\in\NN}$
such that $X^n\in\cC^1([0,T];\m{cB-alpha})$
and converging to $X$ in $\cC^{\gamma'}([0,T];\m{cB-alpha})$
for all $\gamma'<\gamma$.

It follows from \eqref{eq-4} and \eqref{eq-7} that
the speeds of convergence of 
\begin{eqnarray*}
A^\delta\fS(X^n)(s,s';t)&\to& A^\delta\fS(X)(s,s';t)\\
A^\delta\fS(X^n)(s,t;t)&\to& A^\delta\fS(X)(s,t;t)
\end{eqnarray*}
in $\cC^{\gamma'}([0,T];\m{cBdelta})$
are uniform 
in $s\leq s'\leq t$ in $[0,T]$
if $\alpha+\delta<\gamma'$. Hence, we get
that $t\mapsto \fS(X)(s,s';t)$ 
and $t\mapsto \fS(X)(s,t;t)$ are
continuous in $\m{cBdelta}$ for
all $\delta>0$ such that $\alpha+\delta<\gamma$.

\bigskip

\noindent$\circ$ \emph{H\"older continuity of $\fS(X)$.}
By construction,
\begin{equation}
\label{eq-3}
\fS(X)(s,s';t)+\fS(X)(s',s'';t)
=\fS(X)(s,s'';t)
\end{equation}
for all $0\leq s\leq s'\leq s''\leq t$.

Fix $h>0$ and assume that $h=t\ell/2^m$ 
for $\ell\in\event{0,\dotsc,2^m}$ and $m\in\NN$.
Then, by the semi-group property,
\begin{equation*}
\fS^n(0,t;t+h)
=S(h)\sum_{k=0}^{{\lfloor 2^n t/(t+h)\rfloor}-1}
S(t-t_k^n)X(t_k^n,t_{k+1}^n).
\end{equation*}
As $h$ is a dyadic on $[0,t]$, we have
${\lfloor 2^n t/(t+h)\rfloor}=2^{n-m}\ell$.
Passing to the limit
in $\fS^n(ih,(i+1)h;t)$ for $i\in\event{0,\dotsc,\ell-1}$
and using \eqref{eq-3} we obtain that 
\begin{equation}
\label{eq-6}
\fS(X)(0,t;t+h)=S(h)\fS(X)(0,t;t).
\end{equation}
In a similar way, 
\begin{equation}
\label{eq-5}
\fS(X)(t,t+h;t+h)=\fS(X_{t+\cdot})(0,h;h),
\end{equation}
where $X_{t+\cdot}$ is the path $(X_{t+s})_{s\geq 0}$.
Furthermore, the continuity of $h\mapsto \fS(X)(0,t;t+h)$
and $h\mapsto \fS(X)(0,h;h)$ together with 
the continuity of $h\mapsto S(h)x$ for all $x\in\cB$
implies that 
\eqref{eq-6} and \eqref{eq-5} are true even 
if $h$ is not a dyadic point of $[0,t]$.

Now,
\begin{multline*}
A^\delta \fS(X)(0,t+h;t+h)-A^\delta\fS(X)(0,t;t)\\
=A^\delta\fS(X)(0,t;t+h)+A^\delta\fS(t,t+h;t+h)-A^\delta\fS(X)(0,t;t)\\
=(S(h)-\m{Id})A^\delta\fS(X)(0,t;t)+A^\delta\fS(t,t+h;t+h).
\end{multline*}
Using \eqref{eq-4} and \eqref{eq-7}, there exist some constants 
$K_1$ and $K_2$ depending only on $\alpha$, $\beta$, $\delta$ 
and $T$
such that 
\begin{equation*}
\norm[]{}{A^\delta\fS(X)(0,t;t)}
\leq  K_1H(Y)t^{\gamma-\delta-\alpha}
+\norm{}{A^{\delta+\alpha}S(t)Y(0,t)}
\leq K_2H(Y)t^{\gamma-\delta-\alpha}
\end{equation*}
and that, invoking \eqref{eq-5}, 
\begin{equation*}
\norm[]{}{A^\delta\fS(X)(t,t+h;t+h)}
\leq K_2H(Y)h^{\gamma-\delta-\alpha}.
\end{equation*}
Thus, for all $\beta$ in $(0,1]$ such that 
$\alpha+\beta+\delta<\gamma$,
\begin{align*}
\norm{}{(S(h)-\m{Id})A^\delta\fS(X)(0,t;t)}
&\leq K_\beta h^\beta
\norm{}{A^{\delta+\beta}\fS(t,t+h;t+h)}\\
&\leq K_\beta K_2H(Y)h^{\beta}t^{\gamma-\delta-\alpha}
\end{align*}
and hence 
\begin{equation*}
\norm{}{A^\delta \fS(X)(0,t+h;t+h)-A^\delta\fS(X)(0,t;t)}
\leq K_3H(Y)h^{\beta}
\end{equation*}
for all $t\in[0,T-h]$, where $K_3$ is a constant 
that depends only on $T$ and $\alpha$, $\beta$, $\gamma$, $\delta$.
\end{proof}

\section{The non-linear Cauchy problem}\label{nonlin}

We will now define and solve our evolution equation
in the non-linear case:
Let $\cU$ be a separable Banach space, 
and $X$ a $\gamma$-H\"older continuous path with value in $\cU$.
Fix $\delta>0$, $\rho\in\RR$ and let $B$ be a map from $\cB_\delta$
to $\cl(\cU;\cB_\rho)$ (the space of linear bounded operators form
$\cU$ to $\cB_\rho$).

Consider the non-linear Cauchy problem
\begin{equation}
\label{eq-nl-cauchy}
Y(t)=y+\int_0^t B(Y(s))\vd X(s)
-\int_0^t AY(s)\vd s, \qquad y \in \cB_{\delta+\kappa}
\end{equation}
where the solution $Y$ is assumed to be
$\kappa$-H\"older continuous 
with values in $\cB_\delta$. The condition $y \in \cB_{\delta+\kappa}$
on the initial condition is the natural one to ensure that the path
$ t \mapsto S(t)y$ belongs to $\cC^\kappa([0,T],\cB_\delta)$ for any
$T > 0$.
\medskip

The integral with respect $\vd X$ in the r.h.s. of
eq.(\ref{eq-nl-cauchy}) is a Young integral which must be understood
according to the following proposition:
\begin{proposition}
\label{prop-young}
Let $\cV$ be a Banach space and let $H \in \cC^\alpha([0,T];L(\cU,\cV))$ with $\alpha + \gamma > 1$.
Then the integral
\begin{equation}
\label{eq-young-integral}
F(t) = \int_0^t H(s) \vd X(s)  
\end{equation}
exists as the limit in $\cV$ of the sums 
$$
F^\Pi(t) = \sum_{\{s_i \}} H(s_i) [X({s_{i+1}})-X({s_{i}})]
$$
over the partitions $\Pi = \{s_i, i=0,\dots,n : 0=s_0\le s_1 \dots \le
s_n = t\}$ of the interval $[0,t]$ and as the size of the partition
goes to zero. Moreover there exists a constant $K_{\alpha+\gamma}$
depending only on $\alpha + \gamma$ such that
\begin{multline}
  \label{eq-young-bound}
\norm{cV}{F(t)- F(s) - H(s)(X(t)-X(s))} \le\\ K_{\alpha+\gamma} H_{\gamma,T}(X;\cU)
H_{\alpha,T}(H;L(\cU,\cV)) |t-s|^{\alpha+\gamma}  
\end{multline}
for any $0 \le s \le t \le T$.
\end{proposition}
\begin{proof}
This statement is a particular case of a more general one proved by
Young~\cite{Young}, see for example Lyons~\cite{Lyons,LyonsBook}.  
\end{proof}

We will prove that eq.(\ref{eq-nl-cauchy}) has a solution in the
\emph{mild} sense, i.e. we will prove that, under suitable
assumptions, there exists a solution $Y$ to the equation
\begin{equation}
  \label{eq-mild-formulation}
 Y(t) = S(t)y + \int_0^t S(t-s) B(Y(s)) \vd X(s)  
\end{equation}
which is the formal variation of constant solution to
eq.(\ref{eq-nl-cauchy}). 

Our approach to the construction of the solution is to recast
eq.(\ref{eq-mild-formulation}) as a fixed point problem for the
application
$$
\Gamma : \cC^{\kappa}([0,T],\cB_\delta) \to \cC^{\kappa}([0,T],\cB_\delta)
$$
defined as
\begin{equation}
  \label{eq-fixed-point-map}
 \Gamma(Y)(t) = S(t)y + \int_0^t S(t-s) B(Y(s)) \vd X(s)  
\end{equation}

First we will show that $\Gamma$ maps a closed ball of $
\cC^{\kappa}([0,T],\cB_\delta)$ into itself, 
and then, assuming a kind of
Lipschitz condition for the operator $B$, we will show that $\Gamma$ is
a contraction in this ball for a small time $\tau$ and obtain a unique
solution up to $\tau$. 
 
\begin{hypothesis}\label{hypc}
Assume that there exists a function $C$ from 
$\m{cBdelta}$ with values in $L(\m{cBdelta}\otimes\cU,\m{cBrho})$
(where $\m{cBdelta}\otimes\cU$ is the tensor product endowed with the
tensor norm, that is a norm $\|\cdot\|_{\m{cBdelta}\otimes\cU}$
such that $\|x\otimes y\|\leq \|x\|_{\m{cBdelta}}\|y\|_{\cU}$
for all $(x,y)\in\m{cBdelta}\times \cU$:
see for example \cite{LyonsBook}, p.~144 for different possibilites 
for such a norm)
such that 
\begin{equation*}
B(y')x-B(y)x=\int_0^1 C(y+\tau(y'-y))(y'-y)\otimes x\vd \tau.
\end{equation*}
for all $y,y'\in \m{cBdelta}$ and all $x\in\cU$.
Also assume that 
$$
M_B(r) = \sup_{\norm{cBdelta}{y} \le r } \norm*{L(\m{cBdelta}\otimes \cU,\m{cBrho})}{C(y)} < \infty 
$$
for any $r$  and that 
for some increasing function $M_C(r)>0$, $\epsilon \in (0,1]$ and 
for all $y,y',z\in \m{cBdelta}$ such that $\norm{cBdelta}{y} \le r$,$\norm{cBdelta}{y'} \le r$
\begin{equation}\label{cdtderivb}
\norm{cBrho}{C(y')z\otimes x-C(y)z\otimes x}
\leq M_C(r) \norm{cBdelta}{y'-y}^\epsilon\norm{cBdelta}{z}\norm{cU}{x}.
\end{equation}  
for any $x$ in $\cU$.
\end{hypothesis}

\begin{lemma}
\label{thm-1}
Under the previous assumptions on $X$ and Hypothesis~\ref{hypc}, 
if 
\begin{equation}
\label{eq:conditions-on-reg}
\gamma+\kappa>1
\text{ and }
0<\kappa<\max\event{\gamma+\rho-\delta,1},
\end{equation}
then there exists $R > 0$ and $\tau>0$ such that
the closed subset of $\cC^{\kappa}([0,\tau];\cB_\delta)$ defined by 
\begin{multline*}
\cW_{\tau}(y,R)\\
:=\ens{x\in\cC^{\kappa}([0,\tau];\cB_\delta)}{
x(0) = y; H_{\kappa,\tau}(x)\leq R, \sup_{s \in [0,\tau]}
\norm*{\cB_{\delta+\kappa}}{x(s)} \le R}.
\end{multline*}
is invariant under $\Gamma$.
Moreover $\tau$ can be chosen independently of  $y$ if 
$$
\sup_{r > 0}
M_B(r)  <\infty.
$$
\end{lemma}

\begin{proof}
Take $ Y \in \cW_{\tau}(y,R)$, and
note that
$$
\sup_{s\in[0,\tau]} \norm{cBdelta}{Y(s)} \le \norm{cBdelta}{y} + R
\tau^\kappa = \tilde R.
$$

In order to estimate $\Gamma(Y)$ we will in
two steps. First we write
\begin{equation}
\label{eq-gamma-1}
\Gamma(Y)(t)=S(t)y+\int_0^t S(t-s)\vd Z_s,
\end{equation}
where $Z$ is given by
\begin{equation*}
Z(t):=\int_0^t B(Y(s))\vd X(s).
\end{equation*}
Then, Proposition \ref{prop-young} and the assumptions on $B$ imply that
$Z$ is $\gamma$-H\"older continuous with values
in $\cB_\rho$, and that there exists some constant $K$
depending only on the parameters such that
for all $0\leq s\leq t\leq \tau$,
\begin{multline}
\label{eq-Z-estimate-1}
\norm{cBrho}{Z^n(t)-Z^n(s)-B(Y^n(s))(X(t)-X(s))}\\
\leq K M_B(\tilde R) H_{\gamma,\tau}(X;\cU)
H_{\kappa,\tau}(Y;\cB_\delta)
|t-s|^{\gamma+\kappa}.
\end{multline}

This in turn imply that $\Gamma(Y)$ is  given 
by \eqref{eq-gamma-1}.
Indeed from the assumption that $\gamma+\kappa>1$ it follows easily that 
\begin{equation*}
\int_0^t S(t-s)B(Y(s))\vd X(s)=\int_0^t S(t-s)\vd Z(s),
\end{equation*}
since $\sum_{i=0}^{k-1} (t_{i+1}-t_i)^{\gamma+\kappa}$
converges to $0$ with the mesh of the partition 
$0\leq t_0\leq \dotsb\leq t_k\leq t$.

The estimate~(\ref{eq-Z-estimate-1}) can be elaborated further to give
\begin{multline*}
\norm{cBrho}{Z(t)-Z(s)}
\leq 
H_{\gamma,T}(X;\cU)(
\norm{cU->Brho}{B(y)}+M_B(\tilde R) H_{\kappa,\tau}(Y;\cB_\delta)\tau^{\kappa})|t-s|^\gamma\\
+KM_B(\tilde R)
H_{\gamma,T}(X;\cU)
H_{\kappa,\tau}(Y;\cB_\delta)
|t-s|^{\gamma+\kappa},
\end{multline*}
where we used the following inequality to bound the supremum of $B(Y(s))$:
\begin{equation*}
  \begin{split}
\sup_{s\in[0,\tau]} \norm{cU->Brho}{B(Y(s))} 
& \le \norm{cU->Brho}{B(y)} + \sup_{s\in[0,\tau]} \norm{cU->Brho}{B(Y(s)) - B(Y(0))}   
\\& \le \norm{cU->Brho}{B(y)} + \sup_{s,t\in[0,\tau]} \norm{cU->Brho}{B(Y(s)) - B(Y(t))}   
\\& \le \norm{cU->Brho}{B(y)} + M_B(\tilde R) \sup_{s,t\in[0,\tau]} \norm{Bdelta}{Y(s) - Y(t)}
\\& \le \norm{cU->Brho}{B(y)} + M_B(\tilde R) \tau^{\kappa} \left(\sup_{s,t\in[0,\tau]} \frac{\norm{Bdelta}{Y(s) - Y(t)}}{|t-s|^\kappa}\right)
\\& = \norm{cU->Brho}{B(y)}+M_B(\tilde R) H_{\kappa,\tau}(Y;\cB_\delta)\tau^{\kappa}.
  \end{split}
\end{equation*}

Note that
$$
H_{\kappa,\tau}(S(\cdot)y;\m{cBdelta}) \le C_{\kappa} \norm*{\cB_{\delta+\kappa}}{y}.
$$

Thus, 
it follows from \eqref{eq-8}  that 
for some constant $C$ (that does not depend on $Y$),
we have 
\begin{multline}
\label{eq-11}
H_{\kappa,\tau}(\Gamma(Y);\m{cBdelta})
\leq C_{\kappa} \norm*{\cB_{\delta+\kappa}}{y}
\\+CH_{\gamma,T}(X;\cU)(\norm{cU->Brho}{B(y)}+
(1+K)M_B(\tilde R) H_{\kappa,\tau}(Y;\m{cBdelta})\tau^{\kappa}).
\end{multline}
By similar arguments, and using eq.(\ref{eq-9}), we obtain also
\begin{multline}
\label{eq-11-bis}
\sup_{s\in[0,\tau]}\norm*{\cB_{\delta+\kappa}}{\Gamma(Y)(s)}
\leq \norm*{\cB_{\delta+\kappa}}{y}
\\+CH_{\gamma,T}(X;\cU)(\norm{cU->Brho}{B(y)}+
(1+K)M_B(\tilde R) H_{\kappa,\tau}(Y;\m{cBdelta})\tau^{\kappa}),
\end{multline}
where we assumed  that $\tau \le 1$ to obtain an expression similar to 
\eqref{eq-11}.
Choose now $\theta \in (0,1)$ and set
$$
R = (1-\theta)^{-1} \left[C_{\kappa}
\norm*{\cB_{\delta+\kappa}}{y}+CH_{\gamma,T}(X;\cU)\norm{cU->Brho}{B(y)}
\right],
$$
and $\tau > 0$ such that
$$
M_B(\norm*{\m{cBdelta}}{y} + R \tau^\kappa) \tau^\kappa \le 
\theta\left[(1+K) C H_{\gamma,T}(X;\cU)\right]^{-1}.
$$

With this choice of $\tau,R$ the left-hand side
in \eqref{eq-11} is smaller than $R$.
This implies that $\Gamma(Y)$ belongs 
to $\cW_{\tau}(y,R)$ if $Y$ belongs to $\cW_{\tau}(y,R)$.

Moreover, note that if we set
$$
M_B = \sup_{r \ge 0} M_B(r) < \infty
$$
then the choice
$$
\tau^\kappa = \frac{\theta}{(1+K)M_B C H_{\gamma,T}(X;\cU)}.
$$
is equally good and independent of $y$.

\end{proof}

\begin{theorem}
\label{thm-2}
Under Hypothesis~\ref{hypc} for $B$ and
conditions~(\ref{eq:conditions-on-reg}) for $\delta$ and $\kappa$, there exists a time $T > 0$
up to which a unique solution $Y$
of eq.(\ref{eq-mild-formulation}) in
$\cC^{\kappa}([0,T];\m{cBdelta})$ exists.
 If $M_B < \infty$, $T$ can be chosen arbitrarily large.
Moreover  the map $X\mapsto Y$ is Lipschitz continuous from
$\cC^\gamma([0,T];\cU)$ to $\cC^\kappa([0,T];\m{cBdelta})$.
\end{theorem}

\begin{proof}
Choose $\tilde \tau,R$ according to Lemma~\ref{thm-1} 
in order to have $\Gamma :
\cW_{\tilde \tau}(y,R) \to \cW_{\tilde \tau}(y,R)$. Our aim is to show
that there exists $\tau \in (0, \tilde \tau]$ such that $\Gamma$ is a
strict contraction in $\cW_{ \tau}(y,R) \subseteq \cW_{\tilde \tau}(y,R)$.   

Let $Y$ and $Y'$ be two paths in $\cW_{ \tau}(y,R)$.
Set
$$
Z(t) = \int_0^t B(Y(s)) \vd X(s)\text{ and }
Z'(t) = \int_0^t B(Y'(s)) \vd X(s).
$$
By Proposition \ref{prop-young} we have the estimate
\begin{equation}
  \label{eq:estimate-z-zprime}
  \begin{split}
& \norm{cBrho}{Z(t)-Z'(t)-Z(s)-Z'(s)}
\\ &\quad \le K H_{\gamma,\tau}(X;\cU) \Big(
\sup_{u \in [0,\tau]} \norm*{L(\cU,\m{cBrho})}{B(Y(u))-B(Y'(u))} |t-s|^\gamma 
\\ &\quad+
H_{\kappa\epsilon,\tau}(B(Y(\cdot))-B(Y'(\cdot));L(\cU,\m{cBrho})) |t-s|^{\gamma+\kappa\epsilon}
\Big)      
  \end{split}
\end{equation}
for any $s,t \in [0,\tau]$. 
Then,
\begin{multline}
  \label{eq:bound-b-y-sup}
\sup_{u \in [0,\tau]} \norm*{L(\cU,\m{cBrho})}{B(Y(u))-B(Y'(u))}
\le M_B(\tilde R) \sup_{u \in [0,\tau]} \norm{cBrho}{Y(u)-Y'(u)} 
\\ \le M_B(\tilde R) \sup_{s,t \in [0,\tau]} \norm{cBrho}{Y(s)-Y'(s)-Y(t)+Y'(t)} 
\\ \le M_B(\tilde R) H_{\kappa,\tau}(Y(s)-Y'(s);\m{cBrho}) \tau^{\kappa}.
\end{multline}

Moreover, calling $\Delta = Y-Y'$ we have
\begin{equation}
  \label{eq:computation-b-y-b-yprime}
  \begin{split}
 B(Y'(t))x&-B(Y(t))x- B(Y'(s))x+B(Y(s))x = 
\\ & = \int_0^1 C(Y(t)+r \Delta(t)) \vd r \Delta(t)\otimes x
\\ & \quad -
\int_0^1 C(Y(s)+r \Delta(s)) \vd r \Delta(s)  \otimes x  
\\ & = \int_0^1 [C(Y(t)+r \Delta(t))-C(Y(s)+r \Delta(s))] \vd r
\Delta(t)\otimes x 
\\ & \quad +
\int_0^1 C(Y(s)+r \Delta(s)) \vd r [\Delta(t)-\Delta(s)]\otimes x    
  \end{split}
\end{equation}
which yields the following estimate:
\begin{equation}
  \label{eq:estimate-b-y-b-yprime}
  \begin{split}
& \norm*{L(\cU,\m{cBrho})}{B(Y(t))-B(Y'(t))-B(Y(s))-B(Y'(s))}  
\le
\\ 
&  \quad 3 M_C(\tilde R)
(H_{\kappa,\tau}(Y;\m{cBdelta})+H_{\kappa,\tau}(Y';\m{cBdelta}))^\epsilon
\sup_{s\in[0,\tau]}
\norm{cBdelta}{Y(s)-Y'(s)} |t-s|^{\kappa\epsilon}
\\ & \quad +  M_B(\tilde R) H_{\kappa,\tau}(Y-Y';\m{cBdelta}) |t-s|^{\kappa}
\\ 
&  \quad \le H_{\kappa,\tau}(Y-Y';\m{cBdelta}) (3 M_C(\tilde R) (2R)^\epsilon
\tau^{\kappa} + M_B(\tilde R) \tau^{(1-\epsilon)\kappa}) |t-s|^{\kappa\epsilon}.
  \end{split}
\end{equation}
Thus
\begin{equation}
H_{\gamma,\tau}(Z-Z'; \m{cBrho})  
\leq K'_\tau  
H_{\kappa,\tau}(Y-Y';\m{cBdelta})
\end{equation}
with
$$
K'_\tau = K_{\gamma+\kappa\epsilon} H_{\gamma,T}(X;\cU) (3 M_C(\tilde R) (2R)^\epsilon
\tau^{\kappa} + M_B(\tilde R) \tau^{(1-\epsilon)\kappa}+M_B(\tilde R) \tau^\kappa) 
$$
going to zero as $\tau \to 0$.

Thanks to Theorem~\ref{thm-1}, it follows that 
there exists a constant $K''_\tau$
decreasing to $0$ with $\tau$ such that 
\begin{equation*}
H_{\gamma,\tau}\left(\Gamma(Y)-\Gamma(Y');\m{cBdelta}\right)
\leq K''_\tau  H_{\kappa,\tau}(Y'-Y;\m{cBdelta})
\end{equation*}
and
\begin{equation*}
\sup_{s \in[0,\tau]}\norm*{\cB_{\delta+\kappa}}{\Gamma(Y)(s)-\Gamma(Y')(s)}
\leq K''_\tau  H_{\kappa,\tau}(Y'-Y;\m{cBdelta}).
\end{equation*}
Thus, if $T$ is small enough, the map $Y \mapsto \Gamma(Y)$
is $K''_\tau$-Lipschitz on some given ball 
$\cW_{\tau}(y,R)$ with $K''_\tau<1$. 
The uniqueness of the solution of eq.(\ref{eq-mild-formulation})  on $[0,T]$
for $T \le \tau$ follows
immediately. 
Furthermore, in the case where $M_B < \infty$ the independence of $\tilde \tau$
from the starting point $y$ allows to prove
iteratively the uniqueness for any small time interval 
$[i\tau,(i+1)\tau]$ for
$i\leq \lfloor T/\tau\rfloor$ and thus deduce existence and uniqueness
of the solution on an arbitrary time interval $[0,T]$.


The Lipschitz continuity of the map $X \mapsto Y$, where $Y$ is the
solution of the Cauchy problem, can be proven along the following
lines:
let $Y$ the solution of the evolution equation
driven by $X$ and $Y'$ that driven by another
path $X' \in \cC^\gamma([0,T],\cU)$, and set
$$
Z(t) = \int_0^t B(Y(s)) \vd X(s), \qquad
Z'(t) = \int_0^t B(Y'(s)) \vd X'(s),
$$
We can decompose $Z-Z'$ as
\begin{multline*}
Z(t) - Z'(t)=   \int_0^t [B(Y(s))-B(Y'(s))] \vd X(s)\\
+ \int_0^t B(Y(s)) \vd [X(s)-X'(s)]
= \Delta^1(t) + \Delta^2(t),
\end{multline*}
noting that all the Young integrals are well defined.
Then
\begin{multline*}
Y(t) - Y'(t) = \int_0^t S(t-s) d[Z(s) - Z'(s)]\\  
= \int_0^t S(t-s) \vd
\Delta^1(s) + \int_0^t S(t-s) \vd
\Delta^2(s).
\end{multline*}
Using now the same kind of estimates as in the beginning 
of this poof, we obtain a bound of the form
$$
H_{\gamma,\tau}(Y-Y'; \m{cBdelta}) \le K'''_{\tau}
H_{\gamma,\tau}(Y-Y'; \m{cBdelta}) + C H_{\gamma,\tau}(X-X'; \cU) 
$$
where $K'''_\tau$ goes to zero as $\tau$ goes to zero. So for a small
time interval $\tau$ we obtain
$$
H_{\gamma,\tau}(Y-Y'; \m{cBdelta}) \le C' H_{\gamma,\tau}(X-X'; \cU) 
$$
Next, some simple argument on H\"older norms (see e.g.~\cite{Gubinelli}) can be used to show that
this estimate is true up to the existence time $T$ of the solution:
$$
H_{\gamma,T}(Y-Y'; \m{cBdelta}) \le C' H_{\gamma,T}(X-X'; \cU). 
$$
\end{proof}

\begin{remark}
The results of Theorems~\ref{thm-1} and \ref{thm-2} are still true
when \eqref{eq-nl-cauchy} is replaced by 
\begin{equation*}
Y(t)=y+\int_0^t B(Y(s))\vd X(s)+\int_0^t F(Y(s))\vd s
-\int_0^t AY(s)\vd Y(s),
\end{equation*}
where $X$ belongs to $\cC^{\gamma}([0,T];\cU)$, and
where $F$ is a Lipschitz function from 
$\m{cBdelta}$ to $\m{cBrho}$.
\end{remark}

\section{The stochastic heat equation in one dimension}\label{exheat}
In this section, we focus our attention on the (formal) equation
\begin{equation}\label{heateq}
\partial_t Y(t,x)=
\Delta Y(t,x)+\si\lp Y(t,x) \rp  \vd X_t(x), 
\quad t\in\ot,\ x\in\ou,
\end{equation}
with Dirichlet boundary conditions and a null initial condition and
we try to find some simple assumptions on the noise $X$ and on
the function $\si:\RR\to\RR$ ensuring the existence and uniqueness of
a mild solution to (\ref{heateq}) in the sense of Theorems
\ref{thm-0} and \ref{thm-2}. Of course our abstract setting could be
applied to  some more general equations, but the heat equation in
dimension one is a good example to see if our general result can be
used in a simple case, and we are also able to compare it with the
abundant existing literature on this equation.  

\vspace{0.2cm}

To get started, let us describe the kind of noise we consider.

\subsection{Fractional Brownian noise}\label{defrac}

Consider $\{ e_n;n\ge 1 \}$ the trigonometrical basis of $L^2_0(\ou)$, that is
$$
e_n(x)=\lp\frac{\pi}{2}\rp^{1/2}\sin\lp 2\pi n x \rp,
\quad n\ge 1, x\in\ou,
$$
and a collection $\{\hwn;n\ge 1\}$ of independent fractional Brownian motions with Hurst parameter $H\in(0,1)$, all defined on the same probability space $(\Omega,\cf,P)$. Recall that $\hwn$ is then a centered Gaussian process with covariance 
$$
E\lc \hwn(s) \hwn(t) \rc
=
\frac12 \lp t^{2H}+s^{2H}-|t-s|^{2H}  \rp, \quad s,t\in\ot.
$$
In the sequel, the spaces $\cb_{\alpha}$ is also identified with
the usual Sobolev spaces $\wo^{2\alpha,p}$ for $\alpha>0$ and a fixed
$p\ge 1$. When $\zeta \in (0,1)$, $\wo^{\zeta,p}$ can be defined by the completion of the smooth functions on $\ou$ with Dirichlet boundary conditions, with respect to the norm
$$
\| \vp \|_{\zeta,p}^{p}
=
|\vp|_{p}^{p}
+
\int_{\ou^2}\frac{|\vp(x)-\vp(y)|^{p}}{|x-y|^{1+\zeta p}} \vd x \vd y.
$$
Recall that for two conjugate exponents $p,q$ and $\zeta\in(0,1)$, 
$\wo^{-\zeta,q}$ is defined as the dual space, in the distribution sense, of
$\wo^{\zeta,p}$. Recall also that for any $\al,\zeta>0$, the operator
$(-\Delta)^{\al}:\wo^{\zeta,p}\to\wo^{\zeta-2\al,p}$, is
one-to-one. Then, we can identify $\cb_{-\alpha}$ with $\wo^{-2\alpha,p}$ when
$\alpha > 0$. 

\vspace{0.2cm}

With these notations in mind, a Gaussian H\"older continuous process with values in a Sobolev space of negative order can easily be constructed:
\begin{lemma}\label{structnoise}
Let $T,\mu$ be two positive real numbers with $\mu \in (0,1)$, and $\{ q_n;n\ge 1 \}$ a collection of positive numbers such that
\begin{equation}\label{cdtqn}
\snu \lp \frac{q_n}{n^\mu} \rp^2<\infty.
\end{equation}
Consider the random field $X$ defined by the formal series
$$
X(t,x)=\snu q_n e_n(x) \hwn(t),\quad t\in\ot,\ x\in\ou.
$$
Then, for any $\gamma<H$, $\mu < \al< 1$, $\hat p\ge 1$, we have 
$X\in \cC^{\gamma}(\ot,\wo^{-\al,\hat p})$ almost surely.
\end{lemma}

\begin{proof}
Let $\gamma<H$, $\al>\mu$, and $\nu$ be a real number such that $\zeta\equiv 2\nu-\al>0$. Since
$$
(-\Delta)^{-\nu}:\wo^{-\al,\hp}\to\wo^{\zeta,\hp}
$$
is a one-to-one operator, it is enough to see that the process 
$V=(-\Delta)^{-\nu} X$ is almost surely an element of 
$\cC^{\gamma}(\ot,\wo^{\zeta,\hat p})$. However, $V$ is defined by the series
$$
V(t,x)=\snu v_n e_n(x) \hwn(t),\quad t\in\ot,\ x\in\ou,
$$
with
$$
v_n=\frac{c_1 q_n}{n^{2\nu}} \qquad \mbox{and}\qquad
c_1=\frac{1}{(2\pi)^{2\nu}}.
$$
Now, if $t_1,t_2\in\ot$ and $x_1,x_2\in\ou$, we have
\begin{multline}\label{hldest}
E\lc |V(t_2,x_2)-V(t_1,x_1)|^2 \rc
\le
2\Big(
E\lc |V(t_2,x_2)-V(t_1,x_2)|^2 \rc\\
+E\lc |V(t_1,x_2)-V(t_1,x_1)|^2 \rc\Big).
\end{multline}
Furthermore,
\begin{multline*}
E\lc |V(t_2,x_2)-V(t_1,x_2)|^2 \rc
=
E\lc \lp\snu v_n e_n(x_2)\lp \hwn(t_2)-\hwn(t_1)  \rp\rp^2 \rc\\
\le \lp \snu v_n^2 e_n^2(x_2) \rp |t_2-t_1|^{2H}
\le c \lp \snu v_n^2  \rp |t_2-t_1|^{2H}.
\end{multline*}
Now, since $2\nu-\mu>0$, it is easily seen from condition (\ref{cdtqn}) that 
\begin{equation}\label{tmest}
E\lc |V(t_2,x_2)-V(t_1,x_2)|^2 \rc\le c |t_2-t_1|^{2H}.
\end{equation}
On the other hand,
\begin{multline*}
E\lc |V(t_1,x_2)-V(t_1,x_1)|^2 \rc
=
E\lc \lp\snu v_n \lp e_n(x_2)-e_n(x_1)\rp \hwn(t_1)  \rp^2\rc\\
\le c T^{2H} \lp \snu n^{2(2\nu-\mu)}v_n^2   \rp |x_2-x_1|^{2(2\nu-\mu)}.
\end{multline*}
However,
$$
\snu n^{2(2\nu-\mu)}v_n^2
=c\snu \lp\frac{q_n}{n^\mu} \rp^2<\infty,
$$
and thus
\begin{equation}\label{spest}
E\lc |V(t_1,x_2)-V(t_1,x_1)|^2 \rc\le c |x_2-x_1|^{2(2\nu-\mu)}.
\end{equation}
Plugging (\ref{tmest}) and (\ref{spest}) into (\ref{hldest}), we get
$$
E\lc |V(t_2,x_2)-V(t_1,x_1)|^2 \rc\le c
\lp |t_2-t_1|^{2H}+  |x_2-x_1|^{2(2\nu-\mu)}\rp,
$$
and since $V$ is a centered Gaussian process, this yields, by a simple application of Kolmogorov's criterion, that $V$ is almost surely H\"older continuous, with exponent $\gamma$ in time and $\zeta$ in space (recall that $\zeta=2\nu-\al<2\nu-\mu$). It is now easily seen that almost surely,
$V\in \cC^{\gamma}(\ot,\wo^{\zeta,\hat p})$, for any $\hp\ge 1$.
\end{proof}

\subsection{The linear case}

In this section, we just see how to read Theorem \ref{thm-0} in our
fractional Brownian context, in order to compare it with the existing results
(see e.g.~\cite{TTV}).

\begin{proposition}\label{linfbm}
Let $T$ be a positive constant, and $p\ge 1$. Assume that $X$ is a fractional Brownian noise defined as in Lemma \ref{structnoise}, with $\mu< 2H$. Then there exists a unique mild solution to the equation
\begin{equation}\label{eqlinfbm}
\partial_t Y(t,x)=\Delta Y(t,x) \vd t + \dd X_t(x),\quad t\in\ot,\ x\in\ou,
\end{equation}
living in the space $\cC^{\kappa}(\ot,\wo^{\delta,p})$ for all $\delta<2H-\mu$,
$2\kappa<2H-\mu-\delta$.
\end{proposition}

\begin{proof}
This is a trivial consequence of Lemma \ref{structnoise} and Theorem \ref{thm-0}.
\end{proof}

\begin{remark}
In the context of Proposition~\ref{linfbm} the minimal condition under which
(\ref{eqlinfbm}) has a function-valued solution is morally $\snu q_n^2
n^{-4H}<\infty$, which is the same 
necessary and sufficient
condition as the one 
found by Tudor, Tindel and
Viens~\cite{TTV} in order 
to have solutions for
the linear Cauchy problem. This tells us that the path-wise method reaches
optimality in the case of time H\"older regularity  greater than
$1/2$, which is a kind of surprise, since the methods of~\cite{TTV} rely
on the isometric properties  associated with the Gaussian stochastic
integrals which are not exploited in the path-wise approach.
Moreover we can handle a
much more general class of noises. It is also worth noticing that, for any
$H>1/2$, the coefficients $q_n=1$ (that is $\mu=1/2+\ep$ for any
$\ep>0$) are consistent with the assumptions of our proposition, which means
that the white noise in space, considered e.g. by Da Prato-Zabczyk \cite{DPZ}
and Walsh \cite{Wa},  can also be considered in our setting for the heat
equation in dimension 1.
\end{remark}

\subsection{The non-linear Cauchy problem}

In this section, we take up the program of giving a simple enough
example of an equation of the type introduced in eq.~(\ref{heateq}), admitting
a function-valued solution. We get the following Theorem.

\begin{theorem}\label{nonheat}
Let $X$ be a noise defined like in Lemma \ref{structnoise}, and $\si\in \mathcal{C}_b^2(\RR)$. Suppose that $H>1/2$ and $\mu>0$ satisfy
\begin{equation}\label{cdthmu}
\mbox{There exist }\delta, \kappa\in(0,1) \mbox{ such that }
2\kappa<2H-\mu-\delta \mbox{ and }
H+\kappa>1.
\end{equation}
Then, there exists a unique mild solution to (\ref{heateq}) in 
$\cC^\kappa(\ot,\wo^{\delta,p})$, up to an explosion time $T$, for any $p\ge 1/\delta$.
\end{theorem}

\begin{proof}
If (\ref{cdthmu}) holds true, one can always find $\mu<\al<1$ such that
$2\kappa<2H-\al-\delta$ and $H+\kappa>1$. In fact, we will prove that the solution $Y_t$ to (\ref{heateq}) lives in $\wo^{\delta,p}$ for $p$ large enough: according to Lemma \ref{structnoise}, recall that $X\in \cC^{\gamma}(\ot,\wo^{-\al,\hat p})$ for any $\gamma< H$, $\al>\mu$ and $\hp\ge 1$. Thus, going back to Theorem~\ref{thm-2}, we will first check that the operator-valued functional, defined on test functions 
$u,\vp\in \mathcal{C}_0^\infty(\ou)$ by
$$
\lc B(u)\vp \rc(t,x)=\si\lp u(t,x) \rp \, \vp(t,x),
$$
can be extended into an application
$$
B:\wo^{\delta,p}\to\cl\lp \wo^{-\al,\hat p};\wo^{-\al,p}  \rp,
$$
for $p$ large enough, and that this application satisfies
Hypothesis~\ref{hypc}.
This will be done in two steps.

\vspace{0.3cm}

\noindent
{\it Step 1:}
If $\vp\in\wo^{-\al,\hat p}$, $B(u)\vp$ is defined by duality: if $\hq$ is the conjugate of $\hp$, we set
$$
_{\wo^{-\al,\hat p}}\Big( B(u)\vp\,;\, \psi\Big)_{\wo^{\al,\hat q}}
=\,
_{\wo^{-\al,\hat p}}\Big( \vp\,;\, B(u)\psi\Big)_{\wo^{\al,\hat q}}.
$$
Hence, a first condition in order to define $B$ properly is, for any $t\in\ot$ , that
\begin{equation}\label{cdtsob}
\si\lp u(t,.) \rp\in \wo^{\al,\hat q}
\quad\mbox{and}\quad
\wo^{\al,\hat q}\mbox{ is an algebra.}
\end{equation}
Now, if $u(t,.)\in \wo^{\al,\hat q}$ and $\si\in \mathcal{C}_b^2(\RR)$, $\si(u(t,.))$
is still an element of $\wo^{\al,\hat q}$, and thanks to the classical Sobolev imbeddings (see e.g. \cite{Ad}), if we assume that $u(t,.)\in\wo^{\delta,p}$, the condition (\ref{cdtsob}) is induced by
$$
\frac1p <\frac{1}{\hq}+\delta-\al
\quad\mbox{and}\quad
\al\hq>1,
$$
and it is readily checked that this condition is equivalent to
\begin{equation}\label{cdtphp}
\hp <\frac{1}{1-\al}
\quad\mbox{and}\quad
p > \frac{1}{\delta}.
\end{equation}
Furthermore, we have seen that $\hp$ (and hence $\hq$) can be chosen
arbitrarily, and thus, given $\al$ and $\delta$, one can always find $p$ and
$\hp$ satisfying condition~(\ref{cdtphp}). We assume that these
coefficients have been chosen and fixed in the remainder of the proof.

\vspace{0.3cm}

\noindent
{\it Step 2:}
Let us check now Hypothesis \ref{hypc}: 
we use again a duality  argument, and hence, it is enough to show that,
if $u,v\in\wo^{\delta,p}$, and $\psi\in \wo^{\al,\hat q}$, then
\begin{equation*}
\left\| [B(u)-B(v)]\psi \right\|_{\wo^{\al,\hat q}}
=
\int_0^1 \lc D(u+\tau(v-u))\rc(v-u)\otimes\psi \, \vd\tau,
\end{equation*}
with
\begin{equation}\label{weaksob}
\|D(h)  \|_{\cl(\wo^{\delta,p}\otimes\wo^{\al,\hat q};\wo^{\al,\hat q})}
\le
c\lp 1+\| h \|_{\wo^{\delta,p}} \rp.
\end{equation}
However, we have been reduced to a situation where $\wo^{\delta,p}$ is continuously imbedded into $\wo^{\al,\hq}$, and thus, inequality (\ref{weaksob}) is implied by 
\begin{equation}\label{difsob}
\left\| [D(h)]w\otimes\psi \right\|_{\wo^{\al,\hat q}}
\le
c \left\| \psi \right\|_{\wo^{\al,\hat q}}
\lp 1+ \left\| h \right\|_{\wo^{\al,\hat q}} 
 \rp
\left\| w \right\|_{\wo^{\al,\hat q}},
\end{equation}
which is the condition we check now: observe first that, from the definition of $B$, we have, for $x\in\ou$,á
$$
[D(h)w\otimes\psi](x)=\si'(h(x))w(x)\psi(x).
$$
Thus, 
the following decomposition holds true:
$$
\left\| [D(h)]w\otimes\psi \right\|_{\wo^{\al,\hat q}}
\le J_1+J_2+J_3+J_4,
$$
with
\begin{eqnarray*}
J_1&=&\left\| \si'(h) \,w \,\psi \right\|_{L^q}\\
J_2&=&
\int_{\ou^2}\frac{|[\si'(h(x))-\si'(h(y))]w(x)\psi(x)|}
{|x-y|^{1+\al\hq}} \vd x\vd y\\
J_3&=&
\int_{\ou^2}\frac{|(w(x)-w(y))\si'(h(y))\psi(x)|}{|x-y|^{1+\al\hq}} \vd x\vd y\\
J_4&=&
\int_{\ou^2}\frac{|(\psi(x)-\psi(y))\si'(h(y))\psi(y)|}{|x-y|^{1+\al\hq}} \vd x\vd y.
\end{eqnarray*}
We now give an upper bound for the term $J_2$: notice first that, thanks to the fact that $\al\hq>1$, the space $\wo^{\al,\hat q}$ is continuously imbedded in $\mathcal{C}_0(\ou)$. Thus
\begin{equation}\label{psiinf} 
\| \psi\|_\infty\le c \| \psi\|_{\wo^{\al,\hat q}}
\quad \mbox{ and }\quad
\| w\|_\infty\le c \| w\|_{\wo^{\al,\hat q}}.
\end{equation}
Now, plugging (\ref{psiinf}) into the definition of $J_2$, it is readily checked that
$$
J_2\le c \|\si^{(2)}\|_\infty
\left\| h \right\|_{\wo^{\al,\hat q}}
\left\| w \right\|_{\wo^{\al,\hat q}}
\left\| \psi \right\|_{\wo^{\al,\hat q}}.
$$
The terms $J_1$, $J_3$ and $J_4$ can now be handled along the same lines, finishing the proof of (\ref{difsob}). Thus, Hypothesis \ref{hypc} is satisfied in our example. In order to apply Theorem \ref{thm-1} in our situation, it remains to verify (\ref{cdtderivb}). This is left to the reader, since it can be done with the same kind of arguments as above.
\end{proof}

\begin{remark}
Notice that for $H>1/2$, one is still allowed to choose $q_n=1$, that is a white noise in time.
\end{remark}

\begin{remark}
It is also easy to construct an equation of the form (\ref{heateq}) for which the coefficient $B$
is non trivial, but still globally Lipschitz, leading thus to an infinite explosion time: 
just choose a test function $\vp:\ou\in\mathbb{R}$, another function $\si\in C_b^2(\mathbb{R})$, and set 
$$
B(w)=\int_0^1 \si(w(x)) \vp(x) \vd x,
$$
which defines a map from $\wo^{\delta,p}$ to 
$\cl( \wo^{-\al,\hat p};\wo^{-\al,p})$ satisfying the desired global Lipschitz conditions. The proof of this claim is left to the reader.
\end{remark}


\newcommand{\SortNoop}[1]{} \providecommand{\urlbib}[1]{\texttt{\small #1}}

\end{document}